\documentclass[letterpaper,letterpaperjournal]{IEEEtran}

\def\subparagraph{} % because IEEE classes don't define this, but titlesec assumes it's present
\usepackage{titlesec}
\titlespacing*{\section}{0pt}{*1}{*1}
\titlespacing{\subsection}{0pt}{*1}{*1}

\renewcommand{\thesubsubsection}{\arabic{subsubsection}}

\titleformat{\subsubsection}[runin]{\itshape}{\thesubsubsection)}{1em}{}
\titlespacing*{\subsubsection}{\parindent}{0pt}{*1}

\usepackage{todonotes}
\usepackage{cite}
\usepackage{float}
\usepackage{graphicx}
\usepackage{amssymb}
\usepackage[cmex 10]{amsmath}
\usepackage{array}
\usepackage{mdwmath}
\usepackage{balance}

\newtheorem{definition}{Definition}
\newtheorem{lemma}{Lemma}

\newtheorem{obtprob}{Problem}

\usepackage{optidef}

\usepackage{xcolor}

\usepackage{tikz}
\usetikzlibrary{shapes,arrows}
\usetikzlibrary{calc,patterns,decorations.pathmorphing,decorations.markings}

\usepackage{todonotes}
\usepackage{xargs}
\newcommandx{\nader}[2][1=]{\todo[linecolor=orange,backgroundcolor=orange!25,bordercolor=orange,author=Nader,#1]{#2}}
\newcommandx{\ilyasse}[2][1=]{\todo[linecolor=orange,backgroundcolor=orange!25,bordercolor=blue,author=ilyasse,#1]{#2}}

\tikzstyle{block} = [draw, fill=blue!20, rectangle, 
    minimum height=3em, minimum width=6em]
\tikzstyle{sum} = [draw, fill=blue!20, circle, node distance=1cm]
\tikzstyle{input} = [coordinate]
\tikzstyle{output} = [coordinate]
\tikzstyle{pinstyle} = [pin edge={to-,thin,black}]

\title{\LARGE \bf Koopman Operator Approximation  under Negative Imaginary Constraints }

\author{M. A. Mabrok, Ilyasse Aksikas and Nader Meskin  %and Matthew A. Garratt 
\thanks{M. A. Mabrok and Ilyasse Aksikas are with Mathematics Program, Department of Mathematics, Statistics and Physics, College of Arts and Sciences, Qatar University, P.O. box 2713, Doha, Qatar. m.a.mabrok@gmail.com, aksikas@qu.edu.qa. Nader Meskin is with Electircal Engineering Department, College of Engineering, Qatar University, nader.meskin@qu.edu.qa.}

}

\begin{document}

\maketitle

\begin{abstract}

Nonlinear Negative Imaginary (NI) systems arise in various engineering applications, such as controlling flexible structures and air vehicles. However, unlike linear NI systems, their theory is not well-developed. In this paper,  we propose a  data-driven   method for learning a lifted linear NI dynamics that approximates a nonlinear dynamical system using the Koopman theory, which is an operator that captures the evolution of nonlinear systems in a lifted high-dimensional space. The  linear matrix inequality that characterizes the NI property is embedded in the Koopman framework, which results in a non-convex optimization problem. To overcome the numerical challenges of solving a non-convex optimization problem with nonlinear constraints,  the optimization  variables are reformatted in order to convert the optimization problem into  a convex one  with the new variables. We compare our method with local linearization techniques and show that our method can accurately capture the nonlinear dynamics and achieve better control performance. Our method provides a numerically tractable solution for learning the Koopman operator under NI constraints for nonlinear NI systems and opens up new possibilities for applying linear control techniques to nonlinear NI systems without linearization approximations.

\end{abstract}

%dynamics
\section{Introduction}

\begin{figure}
  \centering\includegraphics[width=7cm]{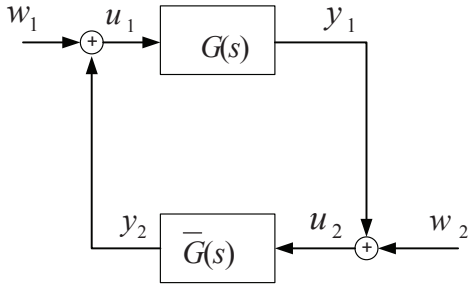}\\
  \caption{If the transfer function matrix of a plant, $G(s)$, is Negative-Imaginary (NI) and the transfer function matrix of a controller, $\bar{G}(s)$, is Strictly Negative-Imaginary (SNI), then the positive-feedback interconnection will be internally stable if and only if the DC gain condition $\lambda_{max}(G(0)\bar{G}(0))<1$ is met.}\label{conn:NI:SNI}
\end{figure}

\IEEEPARstart{N}egative imaginary (NI) systems are systems that have a negative imaginary frequency response when subjected to a positive feedback interconnection \cite{lanzon2008,Mabrok2015TAC,petersen2016negative}. NI systems theory has been widely used for analysis and design of linear-time-invariant (LTI) control systems, especially for flexible structures and air vehicles \cite{lanzon2008,Mabrok2015TAC,vu20202,Preumont2011}. NI systems have several desirable properties, such as robust stability under positive feedback, disspativity  with respect to collocated inputs and outputs, and existence of optimal controllers \cite{petersen2016negative}. A central result in the NI theory is that the positive feedback interconnection between an NI system with  a transfer function matrix $G(s)$ and a strictly NI system with a transfer function matrix $\bar{G}(s)$, as shown in Fig. \ref{conn:NI:SNI}, results in a robustly stable feedback interconnection \cite{lanzon2008,Mabrok2015TAC}. This implies that if a given system is NI, then, it is a great advantage to synthesis a strictly NI controller to guarantee robust stability.

There have been various approaches to address the NI control synthesis problem \cite{mabrok2012stabilization,mabrok2016controller,mabrok2021controller}.  The paper \cite{mabrok2021controller} proposes a synthesis methodology for non-linear systems using a library of controllers parametrized and strictly negative imaginary, optimized through Sequential Quadratic Programming. This methodology can be applied to different non-linear systems. In \cite{mabrok2012stabilization}, a negative imaginary and strict negative imaginary lemma is used to construct a static controller that ensures robust stability against strict negative imaginary uncertainty. In \cite{mabrok2016controller}, a data-driven controller synthesis methodology is proposed that uses measured frequency response data to construct the controller response and transfer function. The data driven linear NI system identification problem was addressed in \cite{mabrok2015system}, where  a modified subspace system identification algorithm that guarantees the negative imaginary property in the identified model by imposing constraints to ensure stability and negative imaginariness was proposed.

Another approach is the sufficient conditions developed in \cite{xiong2016} for designing various types of output feedback controllers, including static, dynamic, and observer-based feedback controllers, based on the above design conditions. The desired controllers can be found using various numerical algorithms. Discrete-time negative imaginary systems were introduced in \cite{ferrante2014foundations,ferrante2017} and the authors investigated the relationship between this concept and discrete-time positive real system theory. Discrete-time negative imaginary systems are defined in terms of a sign condition that must be met in a domain of analyticity of the transfer function. In \cite{liu2017}, similar discrete-time results are explored without the restriction of a symmetric transfer function. 

Negative imaginary systems' theory provides a way to analyze robustness and design robust controllers. In other words,   the theory can be used to design controllers that are able to maintain stability and performance even in the presence of uncertainties in the system, such as unmodeled spillover dynamics or variations in resonant frequencies and damping levels \cite{lanzon2009}.

However, many real-world systems are nonlinear in nature and cannot be adequately modeled by LTI systems. For example, the dynamics of a mass-spring-damper system with a nonlinear spring or a nonlinear damper are nonlinear NI systems \cite{ghallab2022negative}. Nonlinear NI systems pose significant challenges for analysis and control, as the existing NI systems theory does not directly apply to them. Therefore, there is a need to develop a general framework for nonlinear NI systems that can capture their essential features and enable their effective control.

One possible approach to deal with nonlinear systems is to use the Koopman operator, which is a linear operator that describes the evolution of scalar observables of nonlinear systems in an infinite-dimensional Hilbert space \cite{mauroy2020koopman}. The Koopman operator has attracted considerable attention in recent years as a powerful tool for nonlinear dynamics modeling and control \cite{korda2018linear,mauroy2020koopman,bevanda2021koopman}. By using the Koopman operator, one can lift the nonlinear dynamics to a high-dimensional space where they become linear, and then apply linear techniques to analyze and control them. However, learning the Koopman operator from data is generally a challenging task, as it involves solving a non-convex optimization problem with nonlinear constraints \cite{mauroy2020koopman}.

In this paper, we propose a  data-driven method for learning the Koopman operator under NI constraints for nonlinear NI systems. We design  a  data-driven   method for learning a lifted linear NI dynamics that approximates a nonlinear dynamical system using the Koopman theory. To overcome the numerical challenges of solving a non-convex optimization problem with nonlinear constraints, we use a change of variable technique, which allow us to reformulate the problem as a convex optimization that can be solved efficiently. We compare our method with local linearization techniques and show that our method can accurately capture the nonlinear dynamics and achieve better control performance. Our method provides a numerically tractable solution for learning the Koopman operator under NI constraints for nonlinear NI systems and opens up new possibilities for applying linear control techniques to nonlinear systems without linearization approximations.

The rest of the paper is organized as follows. Section II reviews some preliminaries on NI systems theory and Koopman operator theory. Section III presents our main results on learning the Koopman operator under NI constraints for nonlinear NI systems. Section IV illustrates our method on a numerical example; a mass-spring-damper system. Section V concludes the paper and discusses some future directions.

%============================
%============================
%============================

\section{Preliminaries and Notation}
This section introduces the  definitions of  negative imaginary systems and Koopman operator.

\subsection{Negative Imaginary Systems}
\label{sec:NI}

Consider the following LTI system
\begin{align}
\label{eq:xdotn}
&\dot{x}(t) = A x(t)+B u(t), \\
\label{eq:yn} &y(t) = C x(t)+D u(t),
\end{align}%
where $A \in \mathbb{R}^{n \times n},B \in \mathbb{R}^{n \times m},C
\in \mathbb{R}^{m \times n},$ $D \in \mathbb{R}^{m \times m},$ and
with the  square transfer function matrix $G(s)=C(sI-A)^{-1} B+D$.
The transfer function matrix $G(s)$ is said to be strictly proper if
$G(\infty)=D=0$. We will use the notation $
\begin{bmatrix}
\begin{array}{c|c}
A & B \\ \hline C & D
\end{array}
\end{bmatrix}$ to denote the state space realization
\eqref{eq:xdotn}, \eqref{eq:yn}. A NI system is defined as follows:
\begin{definition}\cite{Mabrok2015TAC}\label{Def:NI}
A square transfer function matrix $G(s)$ is NI if  the following
conditions are satisfied:
\begin{enumerate}
\item $G(s)$ has no pole in $\textmd{Re}[s]>0$.
\item For all $\omega >0$ such that $s=j\omega$ is not a pole of $G(s
)$,
\begin{equation}\label{eq:NI:def}
    j\left( G(j\omega )-G(j\omega )^{\ast }\right) \geq 0.
\end{equation}
\item If $s=j\omega _{0}$ with $\omega _{0}>0$ is a pole of $G(s)$, then it is a simple pole and the residue matrix $K=\underset{%
s\longrightarrow j\omega _{0}}{\lim }(s-j\omega _{0})jG(s)$ is
Hermitian and  positive semidefinite.
 \item If $s=0$ is a pole of $G(s)$, then
$\underset{s\longrightarrow 0}{\lim }s^{k}G(s)=0$ for all $k\geq3$
and $\underset{s\longrightarrow 0}{\lim }s^{2}G(s)$ is Hermitian and
positive semidefinite.
\end{enumerate}
\end{definition}
%\end{definition}

\begin{definition}\cite{xiong21010jor}
A square transfer function matrix $G(s)$ is SNI if  the
following conditions are satisfied:
\begin{enumerate}
\item ${G}(s)$ has no pole in $\textmd{Re}[s]\geq0$.
\item For all $\omega >0$, $j\left( {G}(j\omega )-{G}(j\omega )^{\ast }\right) > 0$.
\end{enumerate}
\end{definition}

\subsection{Negative Imaginary lemma}
The negative imaginary lemma is a result that describes NI systems using a pair of LMIs, similar to how the positive-real lemma does \cite{anderson-bk1973,brogliato-bk2007}. This result was presented in \cite{lanzon2007,lanzon2008} and it also covers the case where there are poles on the imaginary axis except at zero \cite{xiong21010jor}.
\begin{lemma}\label{eq:NI-lma}(See \cite{xiong21010jor})
Let $
\begin{bmatrix}
\begin{array}{c|c}
A & B \\ \hline C & D
\end{array}
\end{bmatrix}$
be a minimal state space realization of a transfer function matrix
$G(s)$. Then $G(s)$ is NI if and only if $\textmd{det}(A)\neq0$, $D=D^{T}$
and there exists a real matrix $P>0$ such that
\begin{equation}\label{eq:NI-lmi1}
AP+PA^{*}\leq0, %\ \ B=-APC^{*}.
\end{equation}
 \begin{equation}\label{eq:NI-lmi2}
  B=-APC^{*}.
 \end{equation}
\end{lemma}

 As shown in Fig. \ref{conn:NI:SNI}, there is a positive feedback interconnection between an NI system with transfer function matrix $G(s)$ and an SNI system with transfer function matrix $\bar{G}(s)$. Additionally, let’s assume that the transfer function matrix $G(s)$ has a minimal state space realization,   $
\begin{bmatrix}
\begin{array}{c|c}
A & B \\ \hline C & D
\end{array}
\end{bmatrix},$ and $\bar{G}(s)$ has a minimal state space realization  $
\begin{bmatrix}
\begin{array}{c|c}
\bar{A} & \bar{B} \\ \hline \bar{C} & \bar{D}
\end{array}
\end{bmatrix}.$

\subsection{Nonlinear Negative Imaginary Systems}
The theory of nonlinear NI is not yet well investigated. Preliminary work has been carried out in \cite{ghallab2018extending,mabrok2021dissipativity,ghallab2022negative}. In this subsection, we are presenting the definition of the nonlinear NI.

Consider the following nonlinear dynamical system
\begin{align} %\label{eq:nonlin_nix}
& \dot{x}(t)=f(x(t), u(t)), \\ 
&y(t)=h(x(t)), \label{eq:nonlin_niy}
\end{align}
where $x(t) \in \mathbb{R}^n$, $u(t) \in \mathbb{R}^m$, and $y(t) \in \mathbb{R}^m$ are the system state, control input, and output, respectively.
\begin{definition}\label{def:nonNI}
The nonlinear system  \eqref{eq:nonlin_niy} is said to be  NI if there exists a positive storage function $S_f:\mathbb{R}^n\rightarrow\mathbb{R}$ of a class $C^1$ such that

\begin{equation}\label{lyapd20}
S_f(x(t)) \leq S_f(x(0)) +\int_0^tu(s)^T \dot{y}(s) ds,
\end{equation}
for all $t > 0$. Alternatively, the integral condition \eqref{lyapd20} can be represented by the differential condition  
\begin{equation}\label{lyapd}
    \dot{S_f}(x(t))\leq u(t)^T\dot{y}(t).
\end{equation}
\end{definition}

It is important to note that the  way output measurement is defined in negative imaginary systems is different from passive systems. Negative imaginary systems use position or acceleration measurements, while passive systems use velocity sensors. Negative imaginary systems also allow for transfer functions with a relative degree of up to two, while passive systems permit relative degree zero or one. This is seen in systems that use force actuators and position sensors, like robotics and nano-positioning systems. Stability results for negative feedback interconnections between passive and strictly passive systems do not apply to positive position feedback interconnections between negative imaginary and strictly negative imaginary systems. This was shown in \cite{Petersen2015phy} for physical systems like spring mass damper systems and RLC electrical networks by comparing stability results for feedback interconnections of negative imaginary systems with those for passive systems.

There are several important mechanical systems  that satisfy the  nonlinear negative imaginary property.  For instance, a mass-spring-damper system with positive nonlinear damping coefficient is a nonlinear negative imaginary from force input and displacement output. A more general class of nonlinear negative imaginary  systems is the class of systems that are governed by the following Euler-Lagrangian dynamics,
\begin{equation}\label{eq:EulerL}
M(q) \ddot{q}+C(q, \dot{q}) \dot{q}+G(q)=F.
\end{equation}
It has been shown in \cite{ghallab2022negative} that the dynamical systems governed by Equation \eqref{eq:EulerL} are indeed nonlinear negative imaginary  systems.

\subsection{Koopman theory}
In this section, a brief introduction for the Koopman operator given in \cite{hara2020learning} is presented. The Koopman operator is a mathematical tool used in the analysis of dynamical systems. It was introduced by Bernard Koopman in 1931, and has since become an important tool in the study of nonlinear systems. In particular, in the context of dynamical systems, the Koopman operator provides a powerful tool for studying the evolution of the system over time. By mapping the system's state space to a higher dimensional space, the operator allows for the analysis of nonlinear dynamics through linear methods. This is particularly useful in cases where traditional nonlinear analysis techniques, such as numerical simulations or bifurcation analysis, become computationally expensive or infeasible. One key application of the Koopman operator is in data-driven modeling of complex systems. By using data from the system, the Koopman operator can be estimated and used to approximate the system's behavior. This approach has been particularly successful in fields such as fluid dynamics, where traditional modeling techniques are often impractical due to the complexity of the system.

Consider the following nonlinear difference equation
\begin{align}
\left\{
	\begin{array}{ll}
		x(j+1) &= f\left(x(j),u(j)\right),\\
		y(j)&= g\left(x(j)\right),
	\label{eq.a0}
	\end{array}
\right.
\end{align}
%%%
where $j$ represents the discrete time, $u\in \mathbb{R}^m$ represents the input, $x\in \mathbb{R}^n$ represents the state of the dynamical system, $y\in \mathbb{R}^l$ represents the output, 
and $f(x,u):\mathbb{R}^{n+m} \to \mathbb{R}^{n}$ and $g(x):\mathbb{R}^{n} \to \mathbb{R}^{l}$ are the nonlinear functions. Let $\xi$ to be defined as, 
%%%
$
\xi:=
\left[\begin{array}{c}
x\\
u
\end{array}
\right]
\in \mathbb{R}^{n+m}.
%\label{eq.a1}
$
%%%
Also, define  $\mathbb{F}$ to  be the nonlinear operator given by 
%%%
$
\mathbb{F} (\xi) :=
	\left[\begin{array}{c}
	f(x,u)\\
	\Xi (u)
	\end{array}\right],
%\label{eq.a2.1}
$
%%%
where $\Xi$ is the time-shift operator, i.e $\Xi (u(j)):=u(j+1).$
%%%
%\begin{align}
%\Xi (u(j))&:=u(j+1).
%\nonumber
%\end{align}
This implies that the time evolution of $\xi$ can be given as
%%%
$
\xi(j+1)=\mathbb{F} (\xi(j)).
$
Also, define $\phi_{\mathrm{inf}}(\xi)$ to be  the infinite-dimensional lifting function in the new state space representation, 
%%%%
\begin{align}
\phi_{\mathrm{inf}}(\xi)=\left[\begin{array}{c}
\phi_{1}(\xi)\\
\phi_{2}(\xi)\\
\vdots
\end{array}
\right].
%\label{eq.a4}
%\nonumber
\end{align}
%%%
Now,  we can  define  the Koopman operator $\mathcal{K}$ as
%%%
\begin{align}
\mathcal{K}(\phi_{\mathrm{inf}}(\xi)):=\phi_{\mathrm{inf}}(\mathbb{F} (\xi)),
%\label{eq.a3.2}
\nonumber
\end{align}
%%%
where the time evolution of the lifting function $\phi_{\mathrm{inf}}(\xi)$ is given as follows; 
%%%%
\begin{align}
\phi_{\mathrm{inf}}(\xi(j+1))=\mathcal{K}(\phi_{\mathrm{inf}}(\xi(j))).
\label{eq.a3.1}
%\nonumber
\end{align}
%%%
which is a linear mapping  defined on the infinite-dimensional state space.

In order to approximate the infinite-dimensional Koopman operator by 
 a finite-dimensional approximation,  we define 
the $N_{\phi}$-dimensional lifting function $\phi(\xi):\mathbb{R}^{n+m} \to \mathbb{R}^{N_{\phi}}$ as
%%%
\begin{align}
\phi(\xi)=\left[\begin{array}{c}
\phi_{1}(\xi)\\
\vdots\\
\phi_{N_{\phi}}(\xi)
\end{array}
\right]
\in \mathbb{R}^{N_{\phi}}.
%\label{eq.a6}
\nonumber
\end{align}
%%%
Now, define  $U\in \mathbb{R}^{N_{\phi}\times N_{\phi}}$ to be a finite-dimensional matrix that approximates the Koopman operator $\mathcal{K}$, i.e., $U$ minimizes the following norm 
%%%
\begin{align}
\|U\phi(\xi) - \phi(\mathcal{F}(\xi))\|.
\label{eq.a3.3}
\end{align}
%%%

This implies that 
\begin{align}
\phi(\xi(j+1)) \approx U\phi(\xi(j)),
\label{eq.a3}
\end{align}
%%%
which  describes the behavior of $\phi_{\mathrm{inf}}(\xi)$, defined in  (\ref{eq.a3.1}). Similar to the approach in \cite{hara2020learning}, we  consider  a class of the lifting function 
%$\phi(\xi)$, which takes the following form, 
%%%
\begin{align}
\phi(\xi)=\left[\begin{array}{c}
\psi(x)\\
u
\end{array}
\right]
\in \mathbb{R}^{N+m}, 
\label{eq.a7}
\nonumber
\end{align}
%%%
 where $\psi(x):\mathbb{R}^{n} \to \mathbb{R}^{N}$ is $N$-dimensional lifting function given by
%%%
\begin{align}
&\psi(x)=\left[\begin{array}{c}
\psi_{1}(x)\\
\vdots\\
\psi_{N}(x)
\end{array}
\right]
\in \mathbb{R}^{N},
%\nonumber
\end{align}
%%%
and $N+m=N_{\phi}$ holds.
The matrix $U$ can be given as 
%%%
 \begin{align}
U=\left[\begin{array}{cc}
A_d & B_d
\end{array}
\right]
\in \mathbb{R}^{(N)\times(N+m)},
\nonumber
\end{align}
%%%%
where $A_d\in \mathbb{R}^{N\times N}$ and $B_d \in \mathbb{R}^{N\times m}$.
This implies that  (\ref{eq.a3}) can be rewritten as 
%%%
\begin{align}\label{eq.8.5}
\psi(x(j+1))\approx A_d\psi(x(j))+B_du(j).
\end{align}
%%%
In order to formulate the output equation in Koopman operator,  the approximation of the output equation in (\ref{eq.a0}) can be approximated as follows;
\begin{align*}
y(j) \approx C_d\psi(x(j)),
\end{align*}
%%%
where $C_d \in \mathbb{R}^{l\times N}$. The final linear state-space model can be defined as follows;
%%%
\begin{align}
	\left\{
	\begin{array}{ll}
	\psi(j+1)&=A_d\psi(j)+B_du(j),\\
	y(j)&=C_d\psi(j).
	\end{array}
	\right.
	\label{eq.a10}
\end{align}

\section{Formulating the Koopman NI constraints }
\label{sec:subspace}
In this section,  the NI constraints for discrete dynamical systems are formulated and then imposed  into the learning problem for obtaining the approximated Koopman operator given in \eqref{eq.a3}. Consider the following  discrete time system
\begin{align}
&x(j+1) = A_{d} x(j)+B_{d} u(j),\label{eq:x:ds}\\
 &y(j) = C_{d} x(j)+D_{d} u(j) \label{eq:y:ds},
\end{align}
where $A_{d} \in \mathbb{R}^{n \times n},B_{d} \in \mathbb{R}^{n \times m},C_{d}
\in \mathbb{R}^{m \times n}$ and $D_{d} \in \mathbb{R}^{m \times m}.$ The  NI lemma  given in Lemma  \ref{eq:NI-lma}, which is  formulated in  continuous-time and therefore, a bilinear transformation, as  presented in \cite{mabrok2015system}, in the following form,
\begin{align*}\label{eq:bilinearT}
A &=\frac{1}{T}(I+A_{d})^{-1}(A_{d}-I), \ \ \
B =\frac{1}{\sqrt{T}}(I+A_{d})^{-1}B_{d}, \notag\\
C &=\frac{1}{\sqrt{T}}C_{d}(I+A_{d})^{-1},\ \ \ \ 
D = D_{d}-C_{d}(I+A_{d})^{-1}B_{d},
\end{align*}
with $T$ denoting a sampling time is used to transform the conditions  \eqref{eq:NI-lmi1} into  corresponding discrete time conditions.
The LMI \eqref{eq:NI-lmi1} is transformed from a continuous-time to a discrete-time form  as follows:
\begin{align}
&AP+PA^{T} \leq0 \\
\Leftrightarrow &(I+A_{d})^{-1}(A_{d}-1)P+P(A_{d}^{T}-1)(I+A_{d}^{T})^{-1} \leq0, \notag\\
\Leftrightarrow&(A_{d}P-P)(I+A_{d}^{T})+(P+A_{d}P)(A_{d}^{T}-I) \leq0, \notag\\
\Leftrightarrow&A_{d}P+A_{d}PA_{d}^{T}-P-PA_{d}^{T}+PA_{d}^{T}+A_{d}PA_{d}^{T} \notag\\&\;\;\;-P-A_{d}P \leq0,
\notag\\
\Leftrightarrow&A_{d}PA_{d}^{T}-P+A_{d}PA_{d}^{T}-P \leq0, \notag\\
\Leftrightarrow&A_{d}PA_{d}^{T}-P \leq0. \label{eq:NI-lmi1-ds}
\end{align}

Also, 
\begin{align}
B &=-APC^{T}  \notag\\
%\frac{1}{\sqrt{T}}(I+A_{d})^{-1}B_{d} &=-\frac{1}{T}(I+A_{d})^{-1}(A_{d}-I)P%
%\frac{1}{\sqrt{T}}(I+A_{d}^{T})^{-1}C^{d^{T}} \notag\\
\Leftrightarrow B_{d}& =-\frac{1}{T}(A_{d}-I)P(I+A_{d}^{T})^{-1}C_{d}^{T}. \label{eq:NI-lmi2-ds}
\end{align}
\subsection{Koopman operator learning problem}

Similar approach in \cite{hara2020learning}, define the following data matrices, which are generated from $L$ measurements,
\begin{align}
\Omega & :=[u(j),u(j+1),\dots,u(j+L-1)] \in \mathbb{R}^{m \times L}, \\
Y & :=[y(j),y(j+1),\dots,y(j+L-1)] \in \mathbb{R}^{m \times L}, \\
\Theta & :=[\psi(j),\psi(j+1),\dots,\psi(j+L-1)] \in \mathbb{R}^{N \times L}, \\
\Theta_{+} & :=[\psi(j+1),\psi(j+2),\dots,\psi(j+L)] \in \mathbb{R}^{N \times L}
\end{align}

This implies that the Koopman operator learning problem can be formulated as follows:

\begin{obtprob}
\begin{mini}|s|[3]
{A_d,B_d , C_d}{J_1(A_d, B_d)+J_2(C_d)}
{\label{eq:prob1}}{}
\addConstraint{P  }{\ge 0}
\addConstraint{A_{d}PA_{d}^{T}-P}{ \leq0}
\addConstraint{B_{d}}{ =-\frac{1}{T}(A_{d}-I)P(I+A_{d}^{T})^{-1}C_{d}^{T}}
\end{mini}
where
\begin{align}
J_1(A_d, B_d) & :=\left\|\Theta_{+}-\left[\begin{array}{ll}
A_d & B_d
\end{array}\right]\left[\begin{array}{c}
\Theta \\
\Omega
\end{array}\right]\right\|_F^2, \label{eq:prob1_J1}\\
J_2(C_d) & :=\left\|Y_j-C_d \Theta\right\|_F^2. \label{eq:prob1_J2}
\end{align}
\end{obtprob}

The optimization problem given in \eqref{eq:prob1} is   non-convex due to the non-convexity of the constraints. In order to overcome this issue, we  replace   the matrix inequality  \eqref{eq:NI-lmi1-ds} by a strict inequality with a positive $\alpha>0$ such  that the following holds,
\begin{align}
P &\geq \alpha I \label{eq:NI-lmi2-1-ds:m}\\
A_{d}PA_{d}^{T}-P &\leq-\alpha I. \label{eq:NI-lmi2-ds:m}
\end{align}
Now, using  Schur complements,  the LMIs given in \eqref{eq:NI-lmi2-1-ds:m} and \eqref{eq:NI-lmi2-ds:m} can be written  as follows;

\begin{align}
%P-A_{d}PA_{d}^{T} &>\alpha I  \notag\\
\left[
\begin{array}{cc}
P-\alpha I & A_{d}P \\
PA_{d}^{T} & P%
\end{array}%
\right]  \geq 0 \ \ 
\Rightarrow\left[
 \begin{array}{cc}
P-\alpha I & Q \\
Q^{T} & P%
\end{array}%
\right]  &\geq 0, \label{eq:NI-lmi1-ds-st}
\end{align}
where $Q=A_{d}P$.

The LMI \eqref{eq:NI-lmi1-ds-st} is now convex in the variables $Q$ and $P$. Furthermore,  the matrix $A_{d}$ can be recovered from $A_{d}=QP^{-1}$.
However, the overall optimization problem is still non-convex because the variables in the cost function still ${A_d,B_d , C_d}$. To tackle this issue,  the identification of the matrix $C_d$ is considered as a separate optimization problem by minimizing the cost function given in \eqref{eq:prob1_J2}, which has an optimal solution that is given as 
\begin{equation}\label{eq:cmatrix}
    C_d = Y \Theta^{\dagger},
\end{equation}
where ${(\cdot)}^{\dagger}$ donates  the Moore-Penrose pseudoinverse of $\Theta$.

Second, two weighting matrices $W$ and $\hat{W} $ are introduced into the cost function \eqref{eq:prob1_J1} as follows;

\begin{align}
J_1(A_d, B_d) & :=\left\| W\left(\Theta_{+}-\left[\begin{array}{ll}
A_d & B_d
\end{array}\right] \left[\begin{array}{c}
\Theta \\
\Omega
\end{array}\right] \right)\hat{W}\right\|_F^2, \label{eq:prob12_J1}
\end{align}

where  $\hat{W} $ is selected as:
\begin{align}
\hat{W}= &\begin{bmatrix}
              \Theta \\
                \Omega
         \end{bmatrix}^T \times\left(\begin{bmatrix}
      \Theta \\
        \Omega
         \end{bmatrix}\begin{bmatrix}
              \Theta \\
                \Omega
         \end{bmatrix}^T\right)^{\dag}\times\begin{bmatrix}
                     P & 0 \\
                     0 & I \\
                   \end{bmatrix}.
 \label{eq:W:hat}
\end{align}

It is assumed that there exists  sufficient data to be  collected, which allows us to assume  that $\begin{bmatrix}
              \Theta \\
               \Omega \\
         \end{bmatrix}$ to be full row  rank. This implies that the cost function in \eqref{eq:prob12_J1} can be written as follows:

\begin{align}
J_1(P, Q, B_d) & :=\left\| W\Theta_{+} \hat{W}- W\left[\begin{array}{ll}
Q & B_d
\end{array}\right]\right\|_F^2. \label{eq:prob13_J1}
\end{align}

Therefore, the optimization problem can now be written as follows; 

\begin{obtprob}
\begin{mini}|s|[3]
{P, Q, B_d}{\left\| W\Theta_{+} \hat{W}- W\left[\begin{array}{ll}
Q & B_d
\end{array}\right]\right\|_F^2.}
{\label{eq:prob2}}{}
\addConstraint{\left[
 \begin{array}{cc}
P-\alpha I & Q \\
Q^{T} & P%
\end{array}%
\right] }{\geq 0}.
\end{mini}
\end{obtprob}

where $\hat{W}$ is defined as in \eqref{eq:W:hat}, and the matrix $W$ is a weighting matrix, which  can be treated as a design parameter  to tune the optimization problem.

The optimization problem in \eqref{eq:prob2} is now  convex in the variables ${P, Q, B_d}$. Using any convex optimization solver, such as CVX, can accurately compute the variables ${P, Q, B_d}$, and therefore, we can compute the matrix $A_d$ from the fact that $Q= A_d P$ as $A_d= QP^{-1}$. Note that the matrix $C_d$ is computed separately in \eqref{eq:cmatrix}.

\section{Example: Mass-Spring-Damper System}
Consider a nonlinear mass-spring-damper system 
\begin{equation}\label{eq:springNN}
\ddot{z}(t)+\beta(z(t), \dot{z}(t)) \dot{z}(t)+K(z(t))=u(t), \quad y(t)=z(t).
\end{equation}
where $\beta(z(t), \dot{z}(t))$ and  $k(z(t))$ donate the friction coefficient and  the  spring stiffness, respectively. 

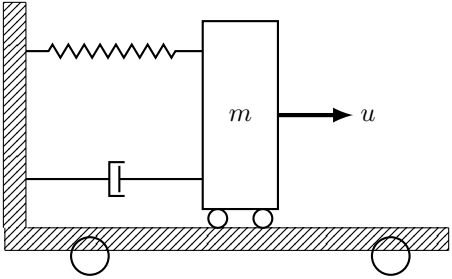
\begin{figure}[H]
\centering
\begin{tikzpicture}
\tikzstyle{spring}=[thick,decorate,decoration={zigzag,pre length=0.3cm,post length=0.3cm,segment length=6}]
\tikzstyle{damper}=[thick,decoration={markings,  
  mark connection node=dmp,
  mark=at position 0.5 with 
  {
    \node (dmp) [thick,inner sep=0pt,transform shape,rotate=-90,minimum width=12pt,minimum height=3pt,draw=none] {};
    \draw [thick] ($(dmp.north east)+(2pt,0)$) -- (dmp.south east) -- (dmp.south west) -- ($(dmp.north west)+(2pt,0)$);
    \draw [thick] ($(dmp.north)+(0,-5pt)$) -- ($(dmp.north)+(0,5pt)$);
  }
}, decorate]
\tikzstyle{ground}=[fill,pattern=north east lines,draw=none,minimum width=0.75cm,minimum height=0.3cm,inner sep=0pt,outer sep=0pt]
\node [style={draw,outer sep=0pt,thick}] (M) [minimum width=1cm, minimum height=2.5cm] {$m$};
\node (ground) [ground,anchor=north,yshift=-0.25cm,minimum width=5.6cm,xshift=-0.03cm] at (M.south) {};
\draw (ground.north east) -- (ground.north west);
\draw (ground.south east) -- (ground.south west);
\draw (ground.north east) -- (ground.south east);
\node (fill) [ground,xshift=-0.15cm,minimum height = 0.3cm, minimum width = 0.3cm] at (ground.west) {};
\draw (fill.north west) -- (fill.south west);
\draw (fill.south west) -- (fill.south east);
\draw [thick] (M.south west) ++ (0.2cm,-0.125cm) circle (0.125cm)  (M.south east) ++ (-0.2cm,-0.125cm) circle (0.125cm);
\draw [thick] (M.south west) ++ (2.5cm,-0.625cm) circle (0.25cm)  (M.south east) ++ (-2.5cm,-0.625cm) circle (0.25cm);
\node (wall) [ground, rotate=-90, minimum width=3cm,yshift=-3cm] {};
\draw (wall.north east) -- (wall.north west);
\draw (wall.north west) -- (wall.south west);
\draw (wall.south west) -- (wall.south east);
\node (y) at (M.east) [xshift = 1.2cm] {$u$};
\draw [spring] (wall.170) -- ($(M.north west)!(wall.170)!(M.south west)$);
\draw [damper] (wall.10) -- ($(M.north west)!(wall.10)!(M.south west)$);
\draw [-latex,ultra thick] (M.east) ++ (0cm,0cm) -- +(1cm,0cm);
\end{tikzpicture}
\caption{Spring-damper-mass system } \label{fig:springmd}
\end{figure}

It can be easily shown that the above system is a nonlinear NI system by choosing  a storage function as follows;

$$
V(z(t), \dot{z}(t))=\frac{1}{2} m \dot{z}^2(t)+\int_0^z(t) k(\xi) d \xi .
$$
The time derivative of $V(z(t), \dot{z}(t))$ gives,
$$
\begin{aligned}
\dot{V}(z(t), \dot{z}(t)) & =m \dot{z}(t) \ddot{z}(t)+k(z(t)) \dot{z}(t)\\&=\dot{z}(t) u(t)-\beta(z(t), \dot{z}(t)) \dot{z}^2(t) 
 \leq \dot{z}(t) u(t),
\end{aligned}
$$
which implies that the spring-damper-mass system  given in \eqref{eq:springNN} is indeed NI system. 

Now,  by defining the states as  ${x}_1(t)=z (t)$ and ${x}_2(t)=\dot{z}(t)$, the spring-damper-mass system,  given in \eqref{eq:springNN}, can be written as follows; 
\begin{align}\label{eq:springSS}
&\left[\begin{array}{l}
\dot{x}_1(t) \\
\dot{x}_2(t)
\end{array}\right]=\left[\begin{array}{c}
x_2(t) \\
\frac{-K(x_1(t))}{m}-\frac{\beta(x_1(t),x_2(t))}{m} x_2(t)+\frac{u(t)}{m}
\end{array}\right] , \\
& y(t)=\left[\begin{array}{ll}
1 & 0
\end{array}\right]\left[\begin{array}{l}
x_1(t) \\
x_2(t)
\end{array}\right] .
\end{align}

Next, the optimization problem in \eqref{eq:prob2} is used to learn a linear system that represents  the nonlinear dynamical system given in \eqref{eq:springSS}. 

The lifting functions  $\phi(x)$ is chosen as a class of  radial basis functions in the following form:
\begin{equation*}
\psi_i(x(j))=\left\|x(j)-r_i\right\|_2^2 \ln \left\|x(j)-r_i\right\|_2,
\end{equation*}
where the center $r_i$ is selected  from the uniform distribution randomly.

The new state space representation in terms of the lifting functions is given as follows, 

\begin{equation}
\psi(x(j))=\left[x_1(j), x_2(j), \psi_1(x(j)),\dots,\psi_{N_{\textmd{rbf}}}(x(j))\right]^{\mathrm{T}} 
\end{equation}

We assumed that  $K(x_1(t))=\left(x_1(t)+x_1^3(t)\right)$ and the nonlinearity in the damper is considered as $\beta(x_1(t),x_2(t)) = x_1^2(t) +x_2^2(t)$. Also, let the number of the lifted states to be 6 states, i.e., $N_{\textmd{rbf}}=6$. The resulting model is compared to the non-constrained Koopman model with the same number of the states, the linearized model around  the initial point  $x_0=[0,0]^{\mathrm{T}}$, and at the linearized model around  the initial point  $x_0=[0.5,0.5]^{\mathrm{T}}$. A random forcing  input signal was used for generating data for the optimization problem.

\begin{figure}[]
    \centering
    \includegraphics[width=0.5\textwidth]{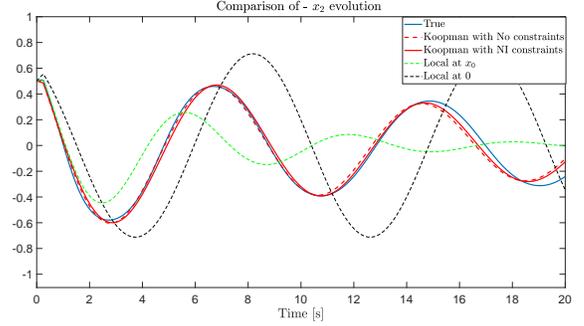}
    \caption{The comparison of the state $x_2$ for different models compared with the true data generated from the nonlinear system}
    \label{fig:state2_6_beta2}
\end{figure}

\begin{figure}[]
    \centering
    \includegraphics[width=0.5\textwidth]{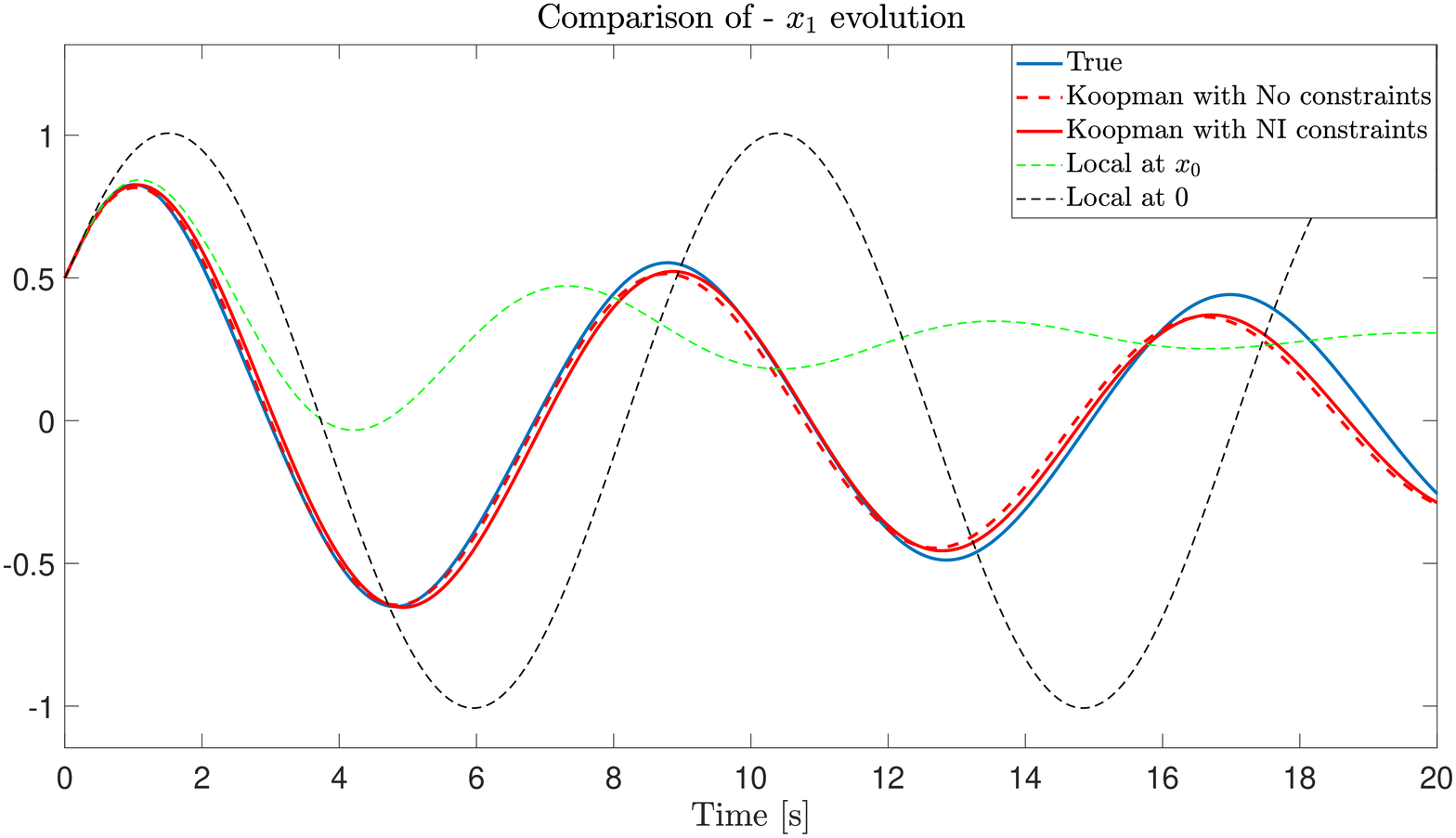}
    \caption{The comparison of the state $x_1$ for different models compared with the true data generated from the nonlinear system}
    \label{fig:state1_6_beta2}
\end{figure}

Figure \ref{fig:state2_6_beta2} presents the state $x_2$ evolutions for different models compared with the true data generated from the nonlinear system. Also, Figure \ref{fig:state1_6_beta2} presents the state $x_1$ evolutions for different models compared with the true data generated from the nonlinear system. As it is shown in Figures \ref{fig:state1_6_beta2} and \ref{fig:state2_6_beta2}, the linear models fail to approximate the nonlinear terms, whereas the Koopman model was able accurately to  approximate the nonlinearity.  

\begin{figure}[]
    \centering
    \includegraphics[width=0.5\textwidth]{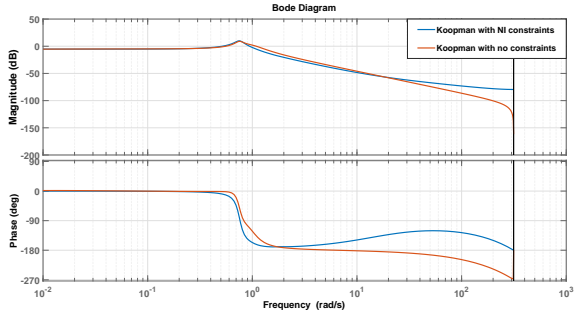}
    \caption{Bode plots for the approximated linear system using Koopman operator. The blue lines represent the NI constrained model developed in this paper. While the orange lines represents the unconstrained model.}
    \label{fig:bode6}
\end{figure}

Figure \ref{fig:bode6} displays the Bode plots corresponding to the approximated linear system employing the Koopman operator. The blue lines depict the model constrained by NI that has been developed in this manuscript, whereas the unconstrained model is represented by the orange lines. It is clear from the phase plot that the NI constrained model satisfy the NI property, i.e., the phase between $(0,-\pi)$, however, the non-constrained model violates the NI property. 

The mean squared error (MSE) between the true values of the states and the generated constrained model is presented in the Table \ref{tab:MSE}. The table shows that the raise in the MSE in the constrained model is very small.   
\begin{table}[]
    \centering
        \caption{The mean squared error for different models.}
    \label{tab:MSE}
    \begin{tabular}{|c|c|c|} \hline
        MSE between:  & $x_1$ & $x_2$ \\  \hline
     True and the constrained model & $0.0078$ & $0.0042$  \\ \hline
      True and the unconstrained model & $0.0030$ & $0.0021$  \\ \hline
       Constrained and unconstrained model & $0.0095$ & $0.0059$  \\ \hline
    \end{tabular}
\end{table}

\begin{figure}[H]
    \centering
    \includegraphics[width=0.5\textwidth]{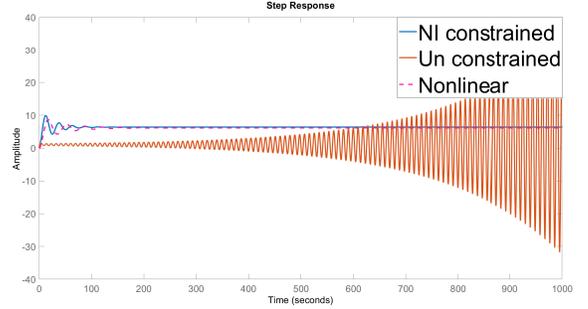}
    \caption{A positive position  feedback controller  connection with the nonlinear model (dashed line), the linear  NI constrained model (blue solid line) and the unconstrained model}
    \label{fig:PPF_control}
\end{figure}

\begin{figure}[H]
    \centering
    \includegraphics[width=0.5\textwidth]{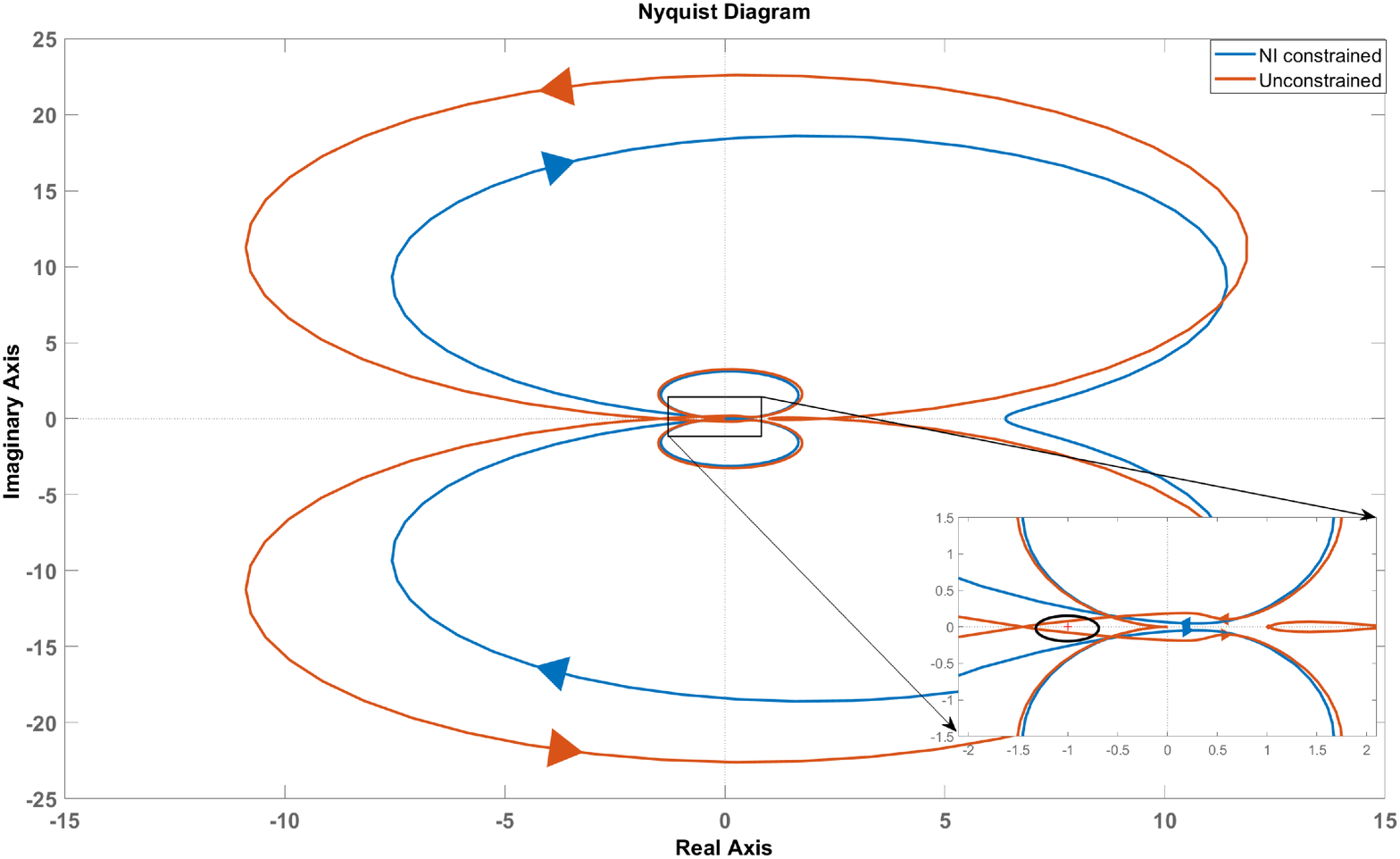}
    \caption{Nyquist plots for the closed loop.}
    \label{fig:Nyq}
\end{figure}

In order to emphasize on  the importance of imposing the NI constraints on the Koopman operator, we consider a feedback control design for the linearized model. Suppose that a strictly NI controller in a form of a positive position feedback controller (PPF) in \eqref{ffp}, which is a common controller for such systems, is designed for the spring-damper-mass. 
\begin{align}\label{ffp}
    C(s) = \frac{K}{s^2 +2\zeta \omega \ s +\omega^2}.
\end{align}
When  the  controller given in   \eqref{ffp} is connected in a positive feedback connection with both, the linear  NI constrained model and the unconstrained model, a stable closed-loop is obtained in the case of the NI constrained model while  in the case of  the unconstrained model, the closed-loop system is unstable. This is shown in the step response of the closed-loop in Figure \ref{fig:PPF_control} and Nyquist plot in Figure \ref{fig:Nyq}. It indeed shows the importance of obtaining a NI linear model approximation of a NI nonlinear system for controller synthesis.

\section{Conclusion }

In conclusion, this paper proposes a data-driven approach to learning a lifted linear Negative Imaginary (NI) dynamics that approximate nonlinear dynamical systems using Koopman theory. The proposed method embeds the linear matrix inequality that characterizes the NI property in the Koopman framework, resulting in a non-convex optimization problem. To overcome the numerical challenges of solving the problem, the optimization variables are reformatted to convert the problem into a convex one. The comparison of the proposed method with local linearization techniques shows that the proposed method can accurately capture the nonlinear dynamics and achieve better control performance. This method provides a numerically tractable solution for learning the Koopman operator under NI constraints for nonlinear NI systems and opens up new possibilities for applying linear control techniques to nonlinear systems without linearization approximations. Overall, the proposed method is a significant contribution to the field of control engineering, as it provides a promising approach for controlling nonlinear NI systems in various engineering applications.

% \bibliographystyle{IEEEtran}
% \bibliography{stability_bib}

\end{document}